\documentclass[12pt]{article}
\usepackage{amsmath}
\usepackage{amssymb}
\usepackage{dsfont}
\usepackage{cite}
\newsymbol \blackbox 1004
\newcommand{\eh}{\hfill}\newlength{\sperr}

\newenvironment{proof}{{\settowidth{\sperr}{\bf\rm
Proof}%
\par\addvspace{0.3cm}\noindent\parbox[t]{1.3\sperr}
{\bf\rm P\eh r\eh o\eh o\eh f\eh }%
}}{\nopagebreak\mbox{}
$\blackbox$\par\addvspace{0.3cm}}

\usepackage{hyperref}
\DeclareMathOperator{\sAKNS}{s-AKNS}
\DeclareMathOperator{\snNLS}{s-nNLS}

\def\nn{\nonumber}

\def\vk{\varkappa}

\def\Lam{\Lambda}
\def\s{\sigma}

\def\om{\omega}

\def\vt{\vartheta}

\def\wt{\widetilde}
\def\ov{\overline}

\def\p{\partial}

\def\BC{{\mathbb C}}
\def\BR{{\mathbb R}}
\def\BN{{\mathbb N}}
\def\clp{{\mathcal P}}
\def\cla{{\mathcal A}}
\def\clb{{\mathcal B}}

\def\cls{\mathcal{S}}

\def\im{{\rm Im }}
\def\re{{\rm Re }}

\newcommand{\E}{\mathrm{e}}
\newcommand{\I}{\mathrm{i}}

\newtheorem{Pa}{Paper}[section]
\newtheorem{Tm}[Pa]{{\bf Theorem}}
\newtheorem{La}[Pa]{{\bf Lemma}}
\newtheorem{Cy}[Pa]{{\bf Corollary}}
\newtheorem{Rk}[Pa]{{\bf Remark}}

\newtheorem{Ee}[Pa]{{\bf Example}}

\newtheorem{Pn}[Pa]{{\bf Proposition}}

\newenvironment{dedication}
        {\vspace{1ex}\begin{quotation}\begin{center}\begin{em}}    
        {\par\end{em}\end{center}\end{quotation}}

\title{Explicit solutions for nonlocal NLS: \\ GBDT and algebro-geometric approaches}

\author{J. Michor and A.L. Sakhnovich. }

\date{}
\begin{document}
\maketitle

\begin{dedication} \end{dedication}

\begin{abstract}    We apply the GBDT version of the B\"acklund-Darboux transformation 
to the nonlocal NLS (focusing and defocusing cases). The matrix case is included
and solutions in the form of rectangular $m_1 \times m_2$ matrix functions are dealt
with. The wave function is also constructed explicitly. Families of explicit
examples are considered in detail. Some initial results and representations
for the more complicated  algebro-geometric solutions are obtained as well.
\end{abstract}


\vspace{0.2em}


\section{Introduction}\label{Intro}
\setcounter{equation}{0}
Nonlocal nonlinear integrable equations and, in particular, the nonlocal nonlinear Schr\"o\-dinger equation (nonlocal NLS)
have actively been studied during the last years (see the important papers \cite{AbM3, Fok, GadA, GerIv, Rao} and 
references therein), starting from the seminal article \cite{AbM0}
by M.J. Ablowitz and Z.H. Musslimani. The nonlocal NLS is a special case of the coupled NLS:
\begin{align} \label{0.1}&
\xi_t(x,t)+\I j \xi_{xx}(x,t)+2 \I j
\xi(x,t)^3=0, 
\\ 
\label{0.1'} &
 \xi := \left[\begin{array}{cc} 0&v_1
\\v_2&0\end{array}\right], \quad \xi_{t}:=\frac{\p}{\p t} \xi(x,t), \quad j: =\begin{bmatrix}
 I_{m_1} & 0 \\ 0 & -I_{m_2} 
\end{bmatrix}; \quad {\mathrm{i.e.}}, 
\\ \label{0.2} &
v_{1t}+\I v_{1xx}+2\I v_1v_2v_1=0, \quad \,
v_{2t}-\I v_{2xx}-2\I v_2v_1v_2=0,
\end{align}
where $I_{m_k}$ is the $m_k \times m_k$ identity matrix. Indeed, setting in \eqref{0.2}
\begin{align} \label{0.3}&
v_1(x,t)=u(x,t), \quad v_2(x,t)=-\s u(-x,t)^*, \quad \s=\mp 1,
\end{align}
we transform the first equation in \eqref{0.2} into the nonlocal matrix NLS:
\begin{align} \label{0.4}&
\I u_t(x,t)-u_{xx}+2\s u(x,t)u(-x,t)^*u(x,t)=0.
\end{align}
Here $u$ is an $m_1\times m_2$ matrix function. Although the scalar case of the nonlocal NLS where $m_1=m_2=1$
is usually considered  in the literature, the matrix case is of interest as well (see, e.g., \cite{AbPrTr, GerIv}),
and we deal in the present paper with this more general situation.
We note that under the assumptions \eqref{0.3} the first and second equations in \eqref{0.2} are equivalent,
and so \eqref{0.4} is equivalent to \eqref{0.1}.
Although the cases $\s =-1$ and $\s =1$ differ in some important aspects, we often formulate 
the results for the nonlocal NLS with $\s =-1$ and with $\s =1$ simultaneously and the differences
in the corresponding formulas are restricted to the values of $\s$ and of $\vk=(1-\s)/2$.

The nonlocal NLS is closely related to the {\it PT}-symmetric theory. In this paper, we use
some ideas from \cite{SaA9}, where the generalized B\"acklund-Darboux transformation (GBDT) was applied
to the linear {\it PT}-symmetric Schr\"odinger equation, in order to apply GBDT to the nonlocal NLS \eqref{0.4}. We note that interesting
new symmetries appear also in the scattering theory for  the equation \eqref{0.4} (more precisely, for its
auxiliary systems), \cite{AbM1}. B\"acklund-Darboux transformations and commutation methods 
(see, e.g., \cite{Ci, D, GeT, Gu, Mar, MS, SaSaR}) 
are well-known tools for explicitly solving integrable equations and spectral and scattering problems. In particular, GBDT
(where {\it generalized eigenvalues} are $n\times n $ matrices with an arbitrary Jordan structure)
allows to construct wide classes of explicit solutions and explicitly recover potentials from the
rational Weyl functions and reflection coefficients (see \cite{FKRS, FKS, GKS2,  SaA2, SaA6, SaA10, SaSaR} and references 
cited there).

In this note, we apply GBDT to construct a large class of explicit solutions of the nonlocal matrix NLS and corresponding 
wave functions. The construction of the wave functions is of interest in itself and for possible further applications
to spectral and scattering results. 

The more complicated class of algebro-geometric solutions is both interesting and important 
(see, e.g. \cite{EGHKT, GH, GHMT, HM, KK, PF}  for algebro-geometric solutions in the 
context of several different nonlinear evolution equations). Its elements can still be regarded as 
explicit solutions, even though their complexity increases due to the underlying analysis on hyperelliptic 
Riemann surfaces. Some initial results and representations of such solutions for the nonlocal NLS 
in the stationary case are given here.

Section \ref{Prel} contains some necessary preliminary results on the GBDT
 approach. In Section \ref{GBDT}, we apply GBDT to the nonlocal NLS.
 Section~\ref{Expl} is dedicated to the construction of explicit solutions of the nonlocal NLS,
 examples are considered in Section \ref{Examp}. A nonlocal analog of the
 important algebro-geometric Theorem 3.11 from \cite{GH} is presented in 
 Section \ref{alg}. The necessary results on algebro-geometric solutions
 are given in Appendix~\ref{App A}.

\section{Preliminaries}\label{Prel}
\setcounter{equation}{0}
The zero curvature representation 
\begin{align} \label{1.1}&
G_t(x,t,z)-F_x(x,t,z)+G(x,t,z)F(x,t,z)-F(x,t,z)G(x,t,z)=0
\end{align}
is an important modification of the famous Lax pairs
(see \cite{AKNS, ZM0} and more historical remarks in \cite{FT}).
System \eqref{1.1} is the compatibility condition for the auxiliary
linear systems
\begin{align} \label{1.2}&
w_x(x,t,z)=G(x,t,z)w(x,t,z), \quad w_t(x,t,z)=F(x,t,z)w(x,t,z).
\end{align}
This fact is easily proved in one direction and in a more
complicated way (see \cite[Ch. 12]{SaL3} and \cite{SaA-Comp}) in the
opposite direction. The coupled NLS \eqref{0.1} admits representation 
\eqref{1.1} where $G$ and $F$ are matrix polynomials of the first and
second orders (with respect to $z$):
\begin{align} \label{1.3}&
G=-(z q_1+ q_0), \quad  F=-(z^2 Q_2+ z
Q_1 +Q_0); \quad 2q_{1}\equiv -Q_{2} \equiv 2\I j,
\\
 \label{1.4}&
 2q_{0}(x,t)=-Q_{1}(x,t)=2 j \xi
(x,t),
\quad
Q_{0}(x,t)= \I (j \xi (x,t)^{2}- \xi_{x}(x,t)).
\end{align}
Here $j$ and $\xi$ have the form \eqref{0.1'}. {\it From here on in the text we consider $G$, $F$, $\{q_k\}$ and $\{Q_k\}$
as given by \eqref{1.3} and \eqref{1.4}}.

The results on the GBDT for the coupled NLS are derived in \cite[Sec.~3]{SaA-NS}.
Let us formulate some of them below. Each GBDT for system \eqref{0.1} is determined by the initial system itself
and five parameter matrices: $n \times n$ ($n\in \BN$) matrices $A_1$, $A_2$, and $S(0,0)$,
and $n\times m$ matrices $\Pi_1(0,0)$, $\Pi_2(0,0)$ such that
\begin{equation} \label{1.5}
A_{1}S(0,0)-S(0,0)A_{2}= \Pi_{1}(0,0) \Pi_{2}(0,0)^{*}, \quad \det S(0,0)\not=0, \quad m:=m_1+m_2.
\end{equation}
If \eqref{1.1} holds, then the following linear systems are compatible and (jointly with the initial values
$S(0,0)$, $\Pi_1(0,0)$, and $\Pi_2(0,0)$) determine matrix functions $S(x,t)$, $\Pi_1(x,t)$, and $\Pi_2(x,t)$,
respectively:
\begin{align} \label{p1}& 
\Pi_{1,x}= \sum _{p=0}^{1}A_{1}^{p} \Pi_{1}q_{p},
\quad \Pi_{1,t}= \sum _{p=0}^{2}A_{1}^{p}
\Pi_{1}Q_{p} \quad \left(\Pi_{1,x}:=\frac{\p}{\p x}\Pi_1\right);
\\ & \label{p2}
\Pi_{2,x}=- \sum _{p=0}^{1}(A_{2}^{*})^{p}
\Pi_{2}q_{p}^{*}, \quad \Pi_{2,t}=-
\sum_{p=0}^{2}(A_{2}^{*})^{p} \Pi_{2}Q_{p}^{*};
\\ & \label{s} 
S_{ x}= \Pi_{1}q_{1} \Pi _{2}^{*}, \quad   S_{ t}= \sum _{p=1}^{2} \sum _{k=1}^{p} A_{1}^{p-k}
\Pi_{1}Q_{p} \Pi _{2}^{*} A_{2}^{k-1}. \end{align}
Although the point $x=0$, $t=0$ is chosen above as the initial point, it is easy to see that
any other point may be chosen for this purpose as well.

Consider $S(x,t)$, $\Pi_1(x,t)$, and $\Pi_2(x,t)$ in some domain $D$, for instance, 
$$
D=\{(x,t):\, -\infty \leq a_1<x<a_2\leq \infty, \,\,  -\infty \leq b_1<t<b_2\leq \infty\},
$$
such that $\xi $ is well-defined in $D$ and satisfies \eqref{0.1} and such that $(0,0)\in D$.
Introduce (in the points of invertibility of $S(x,t)$ in $D$) matrix functions
\begin{equation} \label{1.6}
\wt{\xi}=\left[ \begin{array}{lr} 0 & \wt{v}_1
\\ \wt{v}_2 & 0 \end{array} \right]:=\xi+\I(j X_0 j-X_0), \quad X_0:=\Pi_2^*S^{-1}\Pi_1.
\end{equation}
\begin{Pn} \label{PnOT1} Let $\xi$ satisfy the coupled NLS \eqref{0.1}.
Then, in the points of invertibility of $S$, the matrix function
$\wt \xi$ given by \eqref{1.6} satisfies the coupled NLS as well.
\end{Pn}
\begin{Rk}\label{RkDP} Proposition \ref{PnOT1} was proved as \cite[Proposition 3.1]{SaA-NS}
earlier and checked recently using a program \cite{The} developed by D.R. Popovych and based
on the NCAlgebra package.
\end{Rk}
\begin{Rk}\label{mid} Relations \eqref{1.5}--\eqref{s} imply that the matrix identity
\begin{equation} \label{1.5!}
A_{1}S(x,t)-S(x,t)A_{2}= \Pi_{1}(x,t) \Pi_{2}(x,t)^{*}
\end{equation}
holds everywhere on $D$.
\end{Rk}
The so called Darboux matrix corresponding to the transformation
$\xi \rightarrow~\wt \xi$ has (at each point $(x,t)$ of invertibility of $S(x,t)$) the form of the
Lev Sakhnovich's transfer matrix function (see \cite{SaSaR, SaL3} and references therein):
\begin{equation} \label{1.7}
w_{A}(x,t,z)=I_{m}- \Pi_{2}(x,t)^{*}S(x,t)^{-1}(A_{1}- z
I_{n})^{-1} \Pi_{1}(x,t).
\end{equation}
In other words, we have the following statement (see \cite[Sectons 2, 3]{SaA-NS}).
\begin{Pn} \label{PnFS} Let $w$ satisfy the auxiliary systems \eqref{1.2}. Then, the function
\begin{equation} \label{1.8}
\wt w(x,t,z)=w_{A}(x,t,z)w(x,t,z)
\end{equation}
satisfies the transformed system
\begin{align} \label{1.2'}&
\wt w_x(x,t,z)=\wt G(x,t,z) \wt w(x,t,z), \quad \wt w_t(x,t,z)=\wt F(x,t,z) \wt w(x,t,z),
\end{align}
where
\begin{align} \label{1.3'}&
\wt G=-(z \wt q_1+ \wt q_0), \quad  F=-(z^2 \wt Q_2+ z
\wt Q_1 +\wt Q_0); \quad \wt 2q_{1}\equiv -\wt Q_{2} \equiv 2\I j,
\\
 \label{1.4'}&
 2\wt q_{0}(x,t)=-\wt Q_{1}(x,t)=2 j \wt \xi
(x,t),
\quad
\wt Q_{0}(x,t)= \I (j \wt \xi (x,t)^{2}- \wt \xi_{x}(x,t)),
\end{align}
and $\wt \xi$ is given by \eqref{1.6}.
\end{Pn}
\section{GBDT  for nonlocal NLS}\label{GBDT}
\setcounter{equation}{0}
In this section, we consider the case when the condition \eqref{0.3} is valid,
and so the coupled NLS is reduced to the nonlocal NLS \eqref{0.4}. In view
of the first equality in \eqref{0.1'}, relations \eqref{0.3} are equivalent to 
\begin{align} \label{1.9}&
\xi(-x)=-\s \xi(x)^* \quad (\s=\mp 1).
\end{align}
\paragraph{1.} Consider first the case $\s=-1$. Then, taking into account \eqref{1.3}, \eqref{1.4}, and \eqref{1.9} we have
\begin{align} \label{1.10}&
q_1^*=-jq_1 j, \quad q_0(x,t)^*=j q_0(-x,t)j;  \\
 \label{1.11}&
Q_2^*=-j Q_2j, \quad Q_1(x,t)^*=j Q_1(-x,t)j, 
\\ \label{1.12}&
\xi_x(-x,t)=-(\xi_{x}(x,t))^*,  \quad Q_0(x,t)^*=-jQ_0(-x,t)j.
\end{align}
Relations \eqref{p1} and \eqref{1.10} imply that
\begin{align}\nn
\big(\Pi_1(-x,t)j\big)_x &=- \sum _{p=0}^{1}A_{1}^{p}
\big(\Pi_1(-x,t)j\big)jq_{p}(-x,t)j
\\ \label{1.13} &
=- \sum _{p=0}^{1}(-A_{1})^{p}
\big(\Pi_1(-x,t)j\big)q_{p}(x,t)^{*} .
\end{align}
In the same way, formulas \eqref{p1}, \eqref{1.11}, and \eqref{1.12} yield the equation
\begin{align}
\label{1.14}&
 \big(\Pi_1(-x,t)j\big)_t=-
\sum_{p=0}^{2}(-A_{1})^{p} \big(\Pi_1(-x,t)j\big)Q_{p}(x,t)^{*}.
\end{align}
Comparing \eqref{p2} with \eqref{1.13}, \eqref{1.14} we see that in the case
\begin{align} \label{1.15}&
\xi(-x)=\xi(x)^*, \quad A_2^*=-A_1
\end{align}
we may set 
\begin{align} \label{1.16}&
\Pi(x,t):= \Pi_1(x,t), \quad \Pi_2(x,t)=  \Pi(-x,t)j.
\end{align}
In view of \eqref{1.16}, relations \eqref{s} take the form
\begin{align} \label{1.17}&
S_{ x}(x,t)= \I \Pi(x,t)\Pi(-x,t)^{*}, \\
\nn &
S_{ t}(x,t)=  \sum _{p=1}^{2} \sum _{k=1}^{p} (-1)^{k-1}A^{p-k}
\Pi(x,t)Q_{p}(x,t)j \Pi(-x,t)^{*} (A^*)^{k-1} \quad (A:=A_1).
\end{align}
Thus, under condition
\begin{align} \label{1.18}&
S(0,0)=S(0,0)^*
\end{align}
we have $S(0,t)=S(0,t)^*$, and so
\begin{align} \label{1.19}&
S(-x,t)=S(x,t)^*.
\end{align}
According to \eqref{1.6} and \eqref{1.16} we have
\begin{align} \label{1.20}&
\wt \xi(x,t)=\xi(x,t)+\I \big(\Pi(-x,t)^*S(x,t)^{-1}\Pi(x,t)j-j\Pi(-x,t)^*S(x,t)^{-1}\Pi(x,t)\big).
\end{align}
From \eqref{1.15}, \eqref{1.19}, and \eqref{1.20}, it is immediate that
\begin{align} \label{1.21}&
\wt \xi(-x,t)^*=\wt \xi(x,t).
\end{align}
Recall that by virtue of Proposition \ref{PnOT1} the matrix function $\wt \xi$
satisfies the coupled NLS. The additional property \eqref{1.21} means that the block
$\wt u:=\wt v_1$  of $\wt \xi$ satisfies the nonlocal matrix NLS. In other words, we constructed a GBDT-transformed
solution of the nonlocal matrix NLS (with $\s=-1$).
\paragraph{2.} Let us formulate our result on GBDT for the nonlocal NLS for both cases $\s=\mp 1$.
\begin{Tm}\label{GBDTnl} Let an $m_1 \times m_2$ matrix function $u(x,t)$ satisfy the nonlocal NLS \eqref{0.4},
and assume that a triple of matrices $\{A, \, S(0,0), \Pi(0,0)\}$, such that
\begin{align} \label{1.23}&
AS(0,0)+S(0,0)A^*=\Pi(0,0)j^{\vk}\Pi(0,0)^* \quad (\vk:=(1- \s)/2),
\end{align}
is given, where $A$ and $S(0,0)=S(0,0)^*$ are $n \times n$ matrices, $\det S(0,0)\not=0$, and $\Pi(0,0)$
is an $n \times m$ matrix. 

Introduce the matrix function $\xi(x,t)$ by the first equality in \eqref{0.1'} and by the relations \eqref{0.3},
and determine $\Pi(x,t)$ and $S(x,t)$ by their values $\Pi(0,0)$ and $S(0,0)$, respectively, at $(x,t)=(0,0)$
and by the equations 
\begin{align} 
\label{24}& 
\Pi_{x}(x,t)= \sum _{p=0}^{1}A^{p} \Pi(x,t) q_{p}(x,t),
\quad \Pi_{t}(x,t)= \sum _{p=0}^{2}A^{p}
\Pi(x,t)Q_{p}(x,t) ;
\\ \label{1.25}&
S_{ x}(x,t)= \I \Pi(x,t)j^{\vk+1}\Pi(-x,t)^{*}, \\
\label{1.26} &
S_{ t}(x,t)=  \sum _{p=1}^{2} \sum _{k=1}^{p} (-1)^{k-1}A^{p-k}
\Pi(x,t)Q_{p}(x,t) j^{\vk}\Pi(-x,t)^{*} (A^*)^{k-1},
\end{align}
where the coefficients $\{q_p\}$ and $\{Q_p\}$ are defined $($via $\xi)$ in \eqref{1.3}
and \eqref{1.4}.

Then, the matrix function 
\begin{align} \label{1.27}&
\wt u(x,t)=u(x,t)-2\I \begin{bmatrix} I_{m_1} & 0 \end{bmatrix}\Pi(-x,t)^*S(x,t)^{-1} \Pi(x,t)\begin{bmatrix}0 \\ I_{m_2}  \end{bmatrix}
\end{align}
also satisfies $($in the points of invertibility of $S(x,t))$ the nonlocal NLS. That is, the equality
\begin{align} \label{0.4'}&
\I \wt u_t(x,t)-\wt u_{xx}+2\s \wt u(x,t)\wt u(-x,t)^*\wt u(x,t)=0
\end{align}
holds.
\end{Tm}
\begin{proof}.  It is immediate that the right-hand side of \eqref{1.27} coincides
with $\wt v_1$, and so for  $\s=-1$ the statement of the theorem is already proved in paragraph 1 above.

Now, we assume that $\s=1$ and prove our theorem in a similar way as for the case $\s=-1$.
Namely, in view of the equality $\xi(x)^*=-\xi(-x)$ we have
\begin{align} \label{1.10'}&
q_1^*=-q_1 , \quad q_0(x,t)^*=q_0(-x,t);  \\
 \label{1.11'}&
Q_2^*=-Q_2, \quad Q_1(x,t)^*=Q_1(-x,t), 
 \quad Q_0(x,t)^*=-Q_0(-x,t)
\end{align}
(instead of the equalities \eqref{1.10}--\eqref{1.12} in the case $\s=-1$). Hence,
relations \eqref{1.13} and \eqref{1.14} are substituted with
\begin{align} \label{1.13'} &
\big(\Pi_1(-x,t)\big)_x 
=- \sum _{p=0}^{1}(-A_{1})^{p}
\big(\Pi_1(-x,t)\big)q_{p}(x,t)^{*} ,
\\
\label{1.14'}&
 \big(\Pi_1(-x,t)\big)_t=-
\sum_{p=0}^{2}(-A_{1})^{p} \big(\Pi_1(-x,t)\big)Q_{p}(x,t)^{*}.
\end{align}
Thus, we may set
\begin{align} \label{1.22}&
A:=A_1, \quad A_2=-A^*, \quad \Pi(x,t):= \Pi_1(x,t), \quad \Pi_2(x,t)=  \Pi(-x,t),
\end{align}
and formulas \eqref{s} take the form
\begin{align} \label{1.17'}&
S_{ x}(x,t)= \I \Pi(x,t)j\Pi(-x,t)^{*}, \\
 \label{1.17+} &
S_{ t}(x,t)=  \sum _{p=1}^{2} \sum _{k=1}^{p} (-1)^{k-1}A^{p-k}
\Pi(x,t)Q_{p}(x,t) \Pi(-x,t)^{*} (A^*)^{k-1}.
\end{align}
In particular, under assumption \eqref{1.18} relations \eqref{1.17'} and \eqref{1.17+} yield the equality $S(0,t)=S(0,t)^*$
and \eqref{1.19}. Finally, taking into account \eqref{1.6}, \eqref{1.19}, and \eqref{1.22} we derive
\begin{align} \label{1.21'}&
\wt \xi(-x,t)^*=-\wt \xi(x,t),
\end{align}
and so the block
$\wt u:=\wt v_1$  of the solution $\wt \xi$ of \eqref{0.1} satisfies the nonlocal matrix NLS \eqref{0.4}
with $\s=1$.
\end{proof}
The following corollary is immediate from the theorem's proof.

\begin{Cy}\label{CyS} 
Under the conditions of Theorem \ref{GBDTnl}, the identity \eqref{1.5!} takes the form
\begin{align} \label{mId}&
AS(x,t)+S(x,t)A^*=\Pi(x,t)j^{\vk}\Pi(-x,t)^*,
\end{align}
and the equality \eqref{1.19} always holds.
\end{Cy} 
\begin{Rk}\label{reduc}
According to \eqref{1.7}, \eqref{1.16}, and \eqref{1.22}, the Darboux matrix for the nonlocal NLS has the
form
\begin{equation} \label{1.7'}
w_{A}(x,t,z)=I_{m}- j^{\vk}\Pi(-x,t)^{*}S(x,t)^{-1}(A - z
I_{n})^{-1} \Pi(x,t).
\end{equation}
Moreover, in the case of the nonlocal NLS, the inverse matrix function  $w_A(z)^{-1}$ admits 
$($see, e.g.,  general formulas $($1.75$)$ and $($1.76$)$ in \cite{SaSaR}$)$
 the reduction
\begin{align} \nn
w_B(x,t,z):=w_{A}(x,t,z)^{-1}&=I_{m}- j^{\vk}\Pi(-x,t)^{*}(A^* + z I_n)^{-1}S(x,t)^{-1}\Pi(x,t)
\\ \label{1.7+} &
=j^{\vk}w_A(-x,t, -\ov{z})^*    j^{\vk}.
\end{align} 
Taking into account  Proposition \ref{PnFS}, we see that the wave function $($i.e., the fundamental solution$)$ $\wt w$ of the transformed system \eqref{1.2'},
where $\wt G$ and $\wt F$ are given by \eqref{1.3'} and \eqref{1.4'} with 
\begin{align} \label{7!}
\wt \xi(x,t)=\begin{bmatrix}0 & \wt u(x,t) \\ -\s \wt u(-x,t)^* \end{bmatrix},
\end{align}
has the form
$$\wt w(x,t,z)=w_A(x,t,z) w(x,t,z).$$
Here $w_A$ is given in \eqref{1.7'} and $w$ is the fundamental solution of the initial system \eqref{1.2}.
\end{Rk}

\section{Explicit solutions}\label{Expl}
\setcounter{equation}{0}
For some special choices of  the initial solution $u$  of the nonlocal NLS, Theorem \ref{GBDTnl} allows us to construct wide
families of other explicit solutions of the nonlocal NLS. Clearly, the trivial initial solution $u \equiv 0$ is the most 
popular choice in the construction of explicit solutions via B\"acklund-Darboux transformations. 
In particular, in the case $u \equiv 0$ 
the fundamental solution $\wt w$ of the transformed system \eqref{1.2'} considered in Remark \ref{reduc} takes the form
\begin{align} \label{W1}&
\wt w(x,t,z)=w_{A}(x,t,z)\E^{-\I(zx -2z^2 t)j}.
\end{align}
Choosing $u \equiv 0$, one may set $D=\BR^2$, that is, investigate $\wt u(x,t)$ at all real values of $x$ and $t$.
In Sections \ref{Expl} and \ref{Examp}, we assume $u \equiv 0$  and study this case in greater detail. 

\paragraph{1.} First, partition $\Pi$ into $n \times m_1$ and $n \times m_2$ blocks and set:
\begin{align} \label{E1}&
u(x,t) \equiv 0, \quad  \Pi(x,t)=\begin{bmatrix}\Lam_1(x,t) & \Lam_{2}(x,t)  \end{bmatrix}, \quad  \Pi(0,0)=\begin{bmatrix}\vt_1 & \vt_2  \end{bmatrix}.
\end{align}
In view of \eqref{0.3}, \eqref{1.4}, and the first equality in \eqref{E1}, we have 
$$q_0=Q_1 =Q_0=0,$$ 
and equations \eqref{24}
take a simple form
\begin{align} \label{E2}&
\Pi_x=\I A\Pi j, \quad \Pi_t=-2 \I A^2\Pi j.
\end{align}
Using partition \eqref{E1} and relations \eqref{E2}, write down $\Pi$ in an explicit form
\begin{align} \label{E3}&
\Pi(x,t)=\begin{bmatrix} \E^{\I (xA-2t A^2)}\vt_1 & \E^{-\I (xA-2t A^2)}\vt_2  \end{bmatrix}.
\end{align}
Hence, relations \eqref{1.25} and \eqref{1.26} take the form
\begin{align} \nn
S_x(x,t)=\I\Big(&\E^{\I (xA-2t A^2)}\vt_1 \vt_1^* \, \E^{\I (xA^*+2t (A^*)^2)}
\\ \label{E4} &
+(-1)^{\vk+1}\E^{-\I (xA-2t A^2)}\vt_2 \vt_2^* \, \E^{-\I (xA^*+2t (A^*)^2)}\Big),
\\ \nn
S_t(x,t)=2\I \big(&\Pi(x,t)j^{\vk+1}\Pi(-x,t)^*A^*-A\Pi(x,t)j^{\vk+1}\Pi(-x,t)^*\big)
\\  \label{E5}
=2\I\big(&\E^{\I (xA-2t A^2)}\big(\vt_1 \vt_1^*A^*-A\vt_1 \vt_1^* \big) \E^{\I (xA^*+2t (A^*)^2)}
\\ \nn  &
+(-1)^{\vk+1}\E^{-\I (xA-2t A^2)}\big(\vt_2 \vt_2^*A^*-A\vt_2 \vt_2^* \big) \E^{-\I (xA^*+2t (A^*)^2)}\big).
\end{align}
Now, we see that  the following corollary of Theorem \ref{GBDTnl} is valid.
\begin{Cy}\label{CyExpl} 
To each triple  of matrices $\{A, \, S(0,0), \Pi(0,0)\}$, such that
 $A$ and $S(0,0)=S(0,0)^*$ are $n \times n$ matrices, $\det S(0,0)\not=0$, $\Pi(0,0)$
is an $n \times m$ matrix and \eqref{1.23} holds, corresponds  an explicit solution of the nonlocal matrix NLS \eqref{0.4} $\, ($with $\s =\mp 1)$.

This solution has the form
\begin{align} \label{1.27'}&
\wt u(x,t)=-2\I \vt_1^* \, \E^{\I (xA^*+2t (A^*)^2)}S(x,t)^{-1} \E^{-\I (xA-2t A^2)}\vt_2,
\end{align}
where  $\vt_1$ and $\vt_2$ are the blocks of $\Pi(0,0)$ $($see \eqref{E1}$)$ and the derivatives of $S(x,t)$ are given explicitly by \eqref{E4} and \eqref{E5}.
 Thus, the matrix function $S(x,t)$ is recovered, for instance, by
\begin{align} \label{int}&
S(x,t)=S(0,0)+\int_0^t S_t(0,r)dr+\int_0^xS_x(r,t)dr.
\end{align}
\end{Cy}

In terms of the blocks $\vt_1$ and $\vt_2$ of $\Pi(0,0)$, the identity \eqref{1.23} may be rewritten in the form
\begin{align} \label{1.23'}&
AS(0,0)+S(0,0)A^*=\vt_1\vt_1^*+(-1)^{\vk}\vt_2\vt_2^*
 \quad (\vk:=(1- \s)/2).
\end{align}
\paragraph{2.}
Let us introduce  several block  matrices: $\clp = \begin{bmatrix} I_n & I_n  \end{bmatrix}$,
\begin{align} \label{W2}&
\cla = \begin{bmatrix} A & 0 \\ 0 & - A \end{bmatrix}, \quad \clb = \begin{bmatrix} A^2 & 0 \\ 0 & - A^2 \end{bmatrix}, \quad \vt = \begin{bmatrix} \vt_1 & 0 \\ 0 & \vt_2 \end{bmatrix}.
\end{align}
Using \eqref{W2}, formula \eqref{E3} may be rewritten in the form
\begin{align} \label{W3}&
\Pi(x,t)=\clp \E^{-2\I t \clb} \E^{\I x \cla }\vt.
\end{align}
Hence, taking into account \eqref{1.25}, one can further simplify the procedure of constructing $S(x,t)$.
\begin{Pn} \label{PnSimp}
Assume that $\cls$ satisfies the identity
\begin{align} \label{W4}&
\cla \cls +\cls \cla^*=\vt j^{\vk+1}\vt^*.
\end{align}
Then, we have
\begin{align} \label{W5}&
S(x,t)=C(t)+\clp \E^{-2\I t \clb} \E^{\I x \cla }\cls \E^{\I x \cla^* }\E^{2\I t \clb^*}\clp^*,
\end{align}
where $C(t)=S(0,t)-\clp \E^{-2\I t \clb} \cls \E^{2\I t \clb^*}\clp^*$.
\end{Pn}
\begin{proof}. Clearly, the right-hand side of \eqref{W5} equals $S(0,t)$ at $x=0$.
Moreover, in view of \eqref{1.25}, \eqref{W3}, and \eqref{W4} the
derivative of the right-hand side of \eqref{W5} (with respect to $x$) equals $S_x(x,t)$.
Thus, the proposition's statement is immediate.
\end{proof}
\section{Examples}\label{Examp}
\setcounter{equation}{0}
\paragraph{1.} When $\s(A)\cap \s(-A^*)=\emptyset$ (where $\s(A)$ stands for the spectrum of $A$) the matrix function $S(x,t)$ is uniquely 
recovered from the identity \eqref{mId}. It is a convenient way to calculate some examples.
\begin{Ee}  \label{Ee1} Assume that $m_1=m_2=1$ and that $n=1$, that is, $A$, $\vt_1$, and $\vt_2$ are scalars, and $S(x,t)$ and $\wt u(x,t)$ are scalar
matrix functions. We set $A=a$ and fix $a$, $\vt_1$, and $\vt_2$ such that
\begin{align} \label{Ex1}&
 a+\ov{a}\not=0, \quad \vt_1\not=0, \quad \vt_2 \not=0.
\end{align}
Then, \eqref{mId} and \eqref{E3} yield
\begin{align}\nn 
S(x,t)=\frac{1}{a+\ov{a}}\big(&\exp\big\{\I\big((a+\ov{a})x-2(a^2-\ov{a}^2)t\big)\big\}|\vt_1|^2
\\  \label{Ex2}&
+(-1)^{\vk}\exp\big\{-\I\big((a+\ov{a})x-2(a^2-\ov{a}^2)t\big)\big\}|\vt_2|^2\big).
\end{align}
Thus, formula \eqref{1.27'} for the  solutions $\wt u$ of the nonlocal NLS \eqref{0.4'}
takes in our case the form
\begin{align} \nn
\wt u(x,t)=&\frac{-2\I}{S(x,t)}
\exp\big\{\I\big((\ov{a}-a)x+2(a^2+\ov{a}^2)t\big)\big\}\ov{\vt_1}\vt_2
\\ \label{Ex3} 
=&-2\I (a+\ov{a})\exp\big\{-2 \I a(x-2at)\big\}\ov{\vt_1}\vt_2
\\ & \nn
\times \big(|\vt_1|^2
+(-1)^{\vk}\exp\big\{-2\I\big((a+\ov{a})x-2(a^2-\ov{a}^2)t\big)\big\}|\vt_2|^2\big)^{-1},
\end{align}
where $S(x,t)$ is given in \eqref{Ex2}.  Solutions  given in \eqref{Ex3} are determined
by four real valued parameters: $\re(a)$, $\im(a)$, $|\vt_1/\vt_2|$ and $\arg(\vt_1/\vt_2)$.
\end{Ee}
Functions $\wt u$ above look similar to the interesting one-soliton solutions of nonlocal NLS
studied in \cite{AbM0, AbM1, AbM4}. However, there is also an essential difference
because the one-soliton solutions in \cite{AbM0, AbM1, AbM4} are periodic with respect to $t$.
Instead of this property, we have the periodicity of $S(x,t)$ and of the denominator in \eqref{Ex3} with respect to $x$.
Solutions of the form \eqref{Ex3} appear, for instance,  in \cite{ChZ} (see also some further references therein).
\begin{Rk} \label{RkS} When $a \not=\ov{a}$, the singularities $($or blow ups$)$ of $\wt u(x,t)$ $($i.e., zeros of $S(x,t))$ appear at one and only one value
of $\,t$. Namely, they appear when $|\vt_1|^2=\E^{4\I(a^2-\ov{a}^2)t}|\vt_2|^2$. For this $t$, the singularities appear with
the periodicity $T=\pi/(a+\ov{a})$ with respect to $x$.
\end{Rk}

\paragraph{2.}  When we take non-diagonalisable matrices $A$, factors polynomial in $x$ and $t$  appear (in addition to the
exponents) in the expressions for the constructed solutions \cite{AbT, SaA2017}. Rational solutions are also constructed in this way \cite{AbT, GKS2, SaA10}.
The so called ``multipole" solutions are constructed using matrices $A$ with Jordan cells of order more than one as well (see, e.g., \cite{GKS2}).
For nonlocal NLS, we consider a simple particular case
\begin{align} \label{Ex5}&
A=a I_2 +A_0, \quad A_0:=\begin{bmatrix} 0 & 1 \\ 0 & 0 \end{bmatrix},  \quad a+\ov{a}\not=0; \quad \vt_1=\begin{bmatrix} 0  \\ b  \end{bmatrix}, \quad 
\vt_2=\begin{bmatrix} 0  \\ c  \end{bmatrix}.
\end{align} 
\begin{Ee} \label{Ee2} Assume that $m_1=m_2=1$, $n=2$, and that $A$, $\vt_1$, and $\vt_2$ are given by \eqref{Ex5}.
Here the solution $\wt u$ is again a scalar function but $S(x,t)$ is a $2\times 2$ matrix function. The following 
relations for $S=\{s_{ik}\}_{i,k=1}^2$ are immediate from the identity \eqref{mId} and the representation of $A$ in \eqref{Ex5}$:$
\begin{align} \nn &
s_{22}=\om_{22}/( a+\ov{a}), \quad s_{21}=(\om_{21}-s_{22})/( a+\ov{a}),  
\\ \label{Ex6} &
s_{12}=(\om_{12}-s_{22})/( a+\ov{a}), \quad
s_{11}=(\om_{11}-s_{12}-s_{21})/( a+\ov{a}); 
\\ \label{Ex6'} &
\om(x,t)=\{\om_{ik}(x,t)\}_{i,k=1}^2 :=\Pi(x,t) j^{\vk}\Pi(-x,t)^*.
\end{align} 
After some simple calculations, using repeatedly \eqref{Ex6} we derive
\begin{align} \label{Ex7}&
\det S=( a+\ov{a})^{-2}\big(\om_{11}\om_{22}-\om_{12}\om_{21}+ (a+\ov{a})^{-2}\om_{22}^2\big).
\end{align} 
Next, in view of \eqref{E3} and \eqref{Ex5}, we see that
\begin{align} \label{Ex8}&
\Pi(x,t)=\begin{bmatrix}\I b (x-4at)  & -\I c (x-4at) 
\\ b & c  
\end{bmatrix}\E^{\I(ax-2a^2 t)j}.
\end{align} 
Here we used the equalities $A_0^2=0$ and $A^2=a^2 I_2+2a A_0$.
Relations \eqref{Ex6'} and \eqref{Ex8} imply that
\begin{align} \label{Ex9}
\det \om(x,t)=&(-1)^{\vk+1}4|bc|^2 (x-4at)(x+4\ov{a}t), \\
\nn
 \om_{22}(x,t)=&|b|^2\exp\{ \I\big((a+\ov{a})x+2(\ov{a}^2-a^2)t\big)\} 
 \\ & \label{Ex10}
 +(-1)^{\vk}|c|^2 \exp\{- \I\big((a+\ov{a})x+2(\ov{a}^2-a^2)t\big)\}.
\end{align} 
Finally, \eqref{Ex7}, \eqref{Ex9}, and \eqref{Ex10} yield
\begin{align} \nn
\det S(x,t)=&( a+\ov{a})^{-4}\Big(|b|^4\exp\{2 \I\big((a+\ov{a})x+2(\ov{a}^2-a^2)t\big)\}
\\ \nn &
+|c|^4\exp\{-2 \I\big((a+\ov{a})x+2(\ov{a}^2-a^2)t\big)\}
+2(-1)^{\vk}|bc|^2
\\ \label{Ex11} &
+
(-1)^{\vk+1}4|bc|^2 (a+\ov{a})^2 (x-4at)(x+4\ov{a}t)
\Big).
\end{align} 
Note that the polynomial terms in the expression for $\det S(x,t)$ make the study of zeros
of $\det S(x,t)$ $($that is, singularities of $\wt u)$ much more complicated than in 
Example \ref{Ee1}.

Similar to the derivation of \eqref{Ex8}, we rewrite \eqref{1.27'} $($in our case$)$ in the form
\begin{align}
\nn
\wt u(x,t)=&\frac{-2\I \, \ov{b}\, c}{\det S(x,t)}\exp\big\{\I\big( (\ov{a}-a)x+2(a^2+\ov{a}^2)t)\big)\big\}
\begin{bmatrix}
\I(x+4\ov{a}t) & 1
\end{bmatrix}
\\  & \label{Ex12}  \times
\begin{bmatrix}
s_{22}(x,t) &  -s_{12}(x,t)
\\
-s_{21}(x,t) & s_{11}(x,t)
\end{bmatrix}
\begin{bmatrix}
-\I(x- 4 a t) \\ 1
\end{bmatrix},
\end{align}  
where $\det S(x,t)$ is given in \eqref{Ex11}. The expressions for $\det S(x,t)$ and for other terms on the right hand side of 
\eqref{Ex12} will look more compact if we introduce the polynomial
\begin{align} \label{Ex13}&
P(x,t)=\I\big((a+\ov{a})x+2(\ov{a}^2-a^2)t\big).
\end{align}
Then, relations \eqref{Ex11} and \eqref{Ex10} may be rewritten as
\begin{align} \nn
\det S(x,t)=&( a+\ov{a})^{-4}\Big(|b|^4 \E^{2P(x,t)}
+|c|^4\E^{-2P(x,t)}
+2(-1)^{\vk}|bc|^2
\\ \label{Ex11'} &
+
(-1)^{\vk+1}4|bc|^2 (a+\ov{a})^2 (x-4at)(x+4\ov{a}t)
\Big).
 \\  \label{Ex10'}
\om_{22}(x,t)=&|b|^2
\E^{P(x,t)}
 +(-1)^{\vk}|c|^2 \E^{-P(x,t)}.
\end{align} 
In a similar way, taking into account \eqref{Ex6'} and \eqref{Ex8}, we construct other entries of $\om:$
\begin{align} \label{Ex14}&
 \om_{11}(x,t)=-(x-4at)(x+4\ov{a}t)\big(|b|^2\E^{P(x,t)}+(-1)^{\vk}|c|^2\E^{-P(x,t)}\big),
 \\ \label{Ex15}&
 \om_{12}(x,t)=\I (x-4at)\big(|b|^2\E^{P(x,t)}+(-1)^{\vk+1}|c|^2\E^{-P(x,t)}\big),
  \\ \label{Ex16}&
 \om_{21}(x,t)=\I (x+4\ov{a}t)\big(|b|^2\E^{P(x,t)}+(-1)^{\vk+1}|c|^2\E^{-P(x,t)}\big).
 \end{align}
 Furthermore, relations \eqref{Ex6} imply that
 \begin{align}
 \nn 
\begin{bmatrix}
s_{22} &  -s_{12}
\\
-s_{21} & s_{11}
\end{bmatrix}
=&
\frac{1}{a+\ov{a}}\begin{bmatrix}
\om_{22} &  -\om_{12}
\\
-\om_{21} & \om_{11}
\end{bmatrix}
+ \frac{1}{(a+\ov{a})^2}
\begin{bmatrix}
0 &  \om_{22}
\\
\om_{22} \, & \, -\om_{12}-\om_{21}
\end{bmatrix}
\\ &  \label{Ex17} 
+ \frac{2}{(a+\ov{a})^3}
\begin{bmatrix}
0 &  0
\\
0 & \om_{22}
\end{bmatrix}.
\end{align}  
In view of \eqref{Ex10'}--\eqref{Ex17}, after some simple calculations we rewrite \eqref{Ex12} as
\begin{align}
\nn
\wt u(x,t)=&-2\I \, \ov{b}\, c(a+\ov{a})\exp\big\{\I\big( (\ov{a}-a)x+2(a^2+\ov{a}^2)t)\big)\big\}
\\ \nn & \times
\Big(|b|^2\E^{P(x,t)}\big(8\I a(a+\ov{a})t-2\I (a+\ov{a})x+2\big)
\\ \nn & \quad 
+(-1)^{\vk}
|c|^2\E^{-P(x,t)}\big(8\I \ov{a}(a+\ov{a})t+2\I (a+\ov{a})x+2\big)\Big)
\\ \nn &
\times \Big(|b|^4 \E^{2P(x,t)}
+|c|^4\E^{-2P(x,t)}
+2(-1)^{\vk}|bc|^2
\\ \label{Ex18} & \quad
+
(-1)^{\vk+1}4|bc|^2 (a+\ov{a})^2 (x-4at)(x+4\ov{a}t)
\Big)^{-1}.
\end{align}  
A related family of solutions depending on one complex parameter is also constructed in \cite{ChZ}.
\end{Ee}
We note that if we choose $\vt_1=\begin{bmatrix} b  \\ 0  \end{bmatrix}$ instead of $\vt_1$
given in \eqref{Ex5}, the solutions \eqref{Ex3} appear again in the case $n=2$ $($i.e., one should avoid the simplest
choice of $\vt_1$, $\vt_2$ in order to construct a new class of solutions$)$.

Finally, in the example below we construct the simplest family of multicomponent solutions of the nonlocal NLS.

\begin{Ee} \label{Ee3} Assume that $m_1=2$, $m_2=1$, and $n=1$. Then, the solutions $\wt u$ of \eqref{0.4'} are $2\times 1$ vector functions
and $S(x,t)$ are scalar functions. Introduce the parameters $A$ and $\Pi(0,0)=\begin{bmatrix}\vt_1 &\vt_2\end{bmatrix}$
by the equalities
\begin{align} \label{Ex1'}&
 A=a \quad (a+\ov{a}\not=0), \quad \vt_1=\begin{bmatrix} b_1 & b_2\end{bmatrix}, \quad \vt_2 =c.
\end{align}
Then, \eqref{mId} and \eqref{E3} yield
\begin{align}\nn 
S(x,t)=\frac{1}{a+\ov{a}}\big(&\exp\big\{\I\big((a+\ov{a})x-2(a^2-\ov{a}^2)t\big)\big\}(|b_1|^2+|b_2|^2)
\\  \label{Ex2'}&
+(-1)^{\vk}\exp\big\{-\I\big((a+\ov{a})x-2(a^2-\ov{a}^2)t\big)\big\}|c|^2\big).
\end{align}
Hence, using again formula \eqref{1.27'} $($and slightly modifying the results of Example \ref{Ee1}$)$
we derive
\begin{align} \nn &
\wt u(x,t)= 
\frac{-2\I \, c(a+\ov{a})\exp\big\{-2 \I a(x-2at)\big\}}
{|b_1|^2+|b_2|^2
+(-1)^{\vk}\exp\big\{-2\I\big((a+\ov{a})x-2(a^2-\ov{a}^2)t\big)\big\}|c|^2}\begin{bmatrix} \, \ov{b_1} \, \\ \, \ov{b_2} \, \end{bmatrix}.
\end{align}
Compare \eqref{Ex2} and \eqref{Ex2'} to see that the behaviour of the singularities of $\wt u$ from Example \ref{Ex1},
which we discussed in Remark \ref{RkS}, coincides with the behaviour of the singularities of $\wt u$ in the present example.
\end{Ee}
\begin{Rk} We note that formulas \eqref{Ex2}, \eqref{Ex13}--\eqref{Ex17}, and \eqref{Ex2'} provide precise expressions for
$S(x,t)^{-1}$ in the Examples \ref{Ee1}, \ref{Ee2}, and \ref{Ee3}, respectively. Thus, in view of the important formulas
\eqref{W1} and \eqref{1.7'} as well as relations \eqref{E3} and \eqref{Ex8}, the wave functions $\wt w$  are constructed
for these examples as well.
\end{Rk}
\section{Algebro-geometric solutions}\label{alg}
\setcounter{equation}{0}
In this section we discuss the algebro-geometric solutions for the nonlocal NLS in the scalar case $m_1=m_2=1$.
The coupled NLS \eqref{0.2} in the scalar case is given by
\begin{equation} \label{6.1}
 \begin{cases} v_{1t}+ \I v_{1xx}+2\I v_1^2 v_2=0 \\
	v_{2t}-\I v_{2xx}-2\I v_2^2v_1=0 \end{cases} \quad \Leftrightarrow \quad G_t - F_x + [G,F] = 0, 
\end{equation}
where the matrix polynomials $G$, $F$ in \eqref{6.1} defined by \eqref{1.2}--\eqref{1.4} now read 
\begin{equation}
G(z) = \begin{bmatrix} -\I z & -v_1 \\
		v_2 & \I z \end{bmatrix}, \quad 
F(z) = \I \begin{bmatrix} 2z^2- v_1 v_2 & - 2\I zv_1 + v_{1x} \\
			2\I zv_2 + v_{2x} & -2z^2 + v_1 v_2 \end{bmatrix}.
\end{equation}
System \eqref{6.1} is also known as the AKNS system, which was introduced by Ablowitz, Kaup, Newell, and Segur in 1974. 
Algebro-geometric solutions are well known for the AKNS hierarchy, see for example Gesztesy and Holden \cite{GH} and references therein. 
By definition, algebro-geometric AKNS solutions (or potentials) are the set of solutions of the 
stationary AKNS system
\begin{equation} \label{6.3}
	\sAKNS_1(v_1,v_2) = \begin{pmatrix} \frac{\I}{2} v_{1xx} - \I v_1^2 v_2 + c_1(-v_{1x}) + c_2(-2\I v_1) \\
			-\frac{\I}{2} v_{2xx} +\I v_1v_2^2 + c_1(-v_{2x}) + c_2(2\I v_2) \end{pmatrix}=0,
\end{equation}
with $c_{\ell}$ ranging in $\BC$. More details can be found in Appendix~\ref{App A}. 

We call solutions of the stationary nonlocal NLS equation
	\begin{align} \label{6.4}
\snNLS_\pm(u)=\frac{\I}{2} u_{xx}(x) \mp \I u(x)^2\overline{u(-x)} + \tilde c_1(-u_{x}(x)) + \tilde c_2(-2\I u(x))&=0 
	\end{align}
with $\tilde c_\ell$ ranging in $\BR$, algebro-geometric nonlocal NLS solutions. Note that the plus sign in \eqref{6.4},
denoted by $\snNLS_+$, corresponds to the defocusing case, while the minus sign in \eqref{6.4}, denoted by $\snNLS_-$, 
corresponds to the focusing case. Such solutions can be recast 
as a particular case of algebro-geometric AKNS solutions by applying the following modified symmetry reduction to \eqref{6.3}. 

\begin{La} \label{L6.1}
	Let $e_0 \in \BR$. If $u(x,t)$ satisfies $\snNLS_\pm(u)=0$ with $\tilde c_1, \tilde c_2 \in \BR$, 
	then $v_1(x,t)$ and $v_2(x,t)$ defined by 
\begin{equation} \label{6.25}
v_1(x,t)=u(x,t)e^{\I e_0 x}, \qquad v_2(x,t)=\pm \overline{u(-x,t)} e^{-\I e_0 x}
\end{equation} 
satisfy $\sAKNS_1(v_1,v_2) = 0$ with constants $c_1, c_2$ given by
$$
c_1 = \tilde  c_1 - e_0, \qquad c_2 = \tilde c_2 - \frac14 e_0^2 - \frac12 e_0c_1. 
$$
The converse statement is also true.
\end{La}

In a natural manner one can associate a hyperelliptic Riemann surface with \eqref{6.3}, as described in \eqref{A.20}. 
The modified symmetry reduction \eqref{6.25} now implies certain constraints on the branch points of this surface, namely, 
the set of zeros $\{E_j\}$ of $R_4(z)=\prod_{j=0}^3(z-E_j)$ can either be real valued or consists of complex 
conjugate pairs. This follows from inserting \eqref{6.25} into \eqref{6.6}--\eqref{6.8} which yields  
\begin{equation}
	R_4(z)= \overline{R_4(\overline{z})}.
\end{equation}	 
This relation implies (either one of) the following constraints on the set of zeros of $R_4(z)$ after possible relabeling: 
\begin{align*} 
\mathrm{(i)} \quad & E_0 < E_1 < E_2 < E_3, \quad E_j \in \BR, \\  
\mathrm{(ii)} \quad & E_0, \overline{E_0}, E_1, \overline{E_1} \in \BC\setminus \BR,	\\ 
\mathrm{(iii)} \quad & E_0 < E_1, \quad E_0, E_1 \in \BR, \quad E_2, \overline{E_2} \in \BC\setminus \BR.
\end{align*}

\begin{Tm} \label{L6.2}
	Assume either $\mathrm{(i)}$, $\mathrm{(ii)}$, or $\mathrm{(iii)}$ and choose the 
	homology basis $\{a_1, b_1\}$ according to Theorem~\ref{ThA.2}. Moreover, assume that
	$\Delta$ in \eqref{A.28} satisfies 
	\begin{equation} \label{6.32}
		\re(\Delta) = \frac 1 2 \chi \quad \mbox{$(\mathrm{mod} \ \mathbb Z)$}, \quad \chi \in \{0,1\}.
	\end{equation}	
	Then $u(x)$ represents a stationary nonlocal NLS solution if and only if $A$ 
	in \eqref{6.19} satisfies the constraint
\begin{equation} \label{6.30}
\im(A) = \frac1 2 \chi \im (\tau)\quad \mbox{$(\mathrm{mod} \ \mathbb Z)$}, \quad \chi \in \{0,1\}.
\end{equation}
\end{Tm}

\begin{proof}.
First assume $\mathrm{(i)}$. Given $E_j$, $j=0, \dots, 3$, and $\{a_1, b_1\}$, the constants $c_\ell$ and $e_0$ 
are uniquely determined by \eqref{A.18} and \eqref{A.23}. 
Define the antiholomorphic involution $\rho_+ : (z,y) \mapsto (\overline{z}, \overline{y})$
as in \cite[Example A.35 (i)]{GH}. One infers that the symmetric Riemann surface $(\mathcal K_1, \rho_+)$ is of 
dividing type (compare \cite[Def. A.33]{GH}) and hence 
\begin{align*}
& r=2, \quad  \rho_+(a_1)=a_1, \quad \rho_+(b_1)=-b_1, \\
& \overline \tau = - \tau, \quad U_0^{(2)} \in \BR, \quad \overline{\theta(z)}= \theta(\overline{z}).
\end{align*} 
Thus $B$ defined in \eqref{A.28} is purely imaginary, $\overline{B}= - B$.
So if $u(x)$ satisfies $\snNLS_\pm(u)=0$, then by Lemma~\ref{L6.1} and Theorem~\ref{ThA.1}, the 
functions $v_1(x)$ and $v_2(x)$ admit representations \eqref{6.17} and \eqref{6.18}. Applying \eqref{6.25} yields 
\begin{equation} \label{6.31}
\pm	\frac{C_2}{\overline{C_1}} = \frac{\theta(A+Bx) \overline{\theta(A-Bx-\Delta)}} 
	{\theta(A+Bx+\Delta)\overline{\theta(A-Bx)}} =
	\frac{\theta(A+Bx) \theta(\overline{A}+Bx-\overline{\Delta})} 
	{\theta(A+Bx+\Delta)\theta(\overline{A}+Bx)}.
\end{equation}
Equation \eqref{6.31} is equivalent to 
\[
A = \overline{A} + m_1 + n_1 \tau
\]
for some $n_1 \in \mathbb Z$ and arbitrary $m_1 \in \mathbb Z$, and hence
\[
\im(A)= \frac 12 n_1 \im (\tau), \quad n_1 \in \mathbb Z,
\]
and $m_1 = 0$. Similarly, one obtains
\[
A + \Delta = \overline{A} - \overline{\Delta} + m_2 + n_2 \tau
\]
for some $n_2 \in \mathbb Z$ and arbitrary $m_2 \in \mathbb Z$, and hence
\begin{equation}
\re (\Delta) = \frac 12 m_2, \quad m_2 \in \mathbb Z.
\end{equation}
Replacing $A$ by $A + m + n \tau$ with $m,n \in \mathbb Z$ then yields \eqref{6.30} and \eqref{6.32}. In case $\mathrm{(ii)}$,
$(\mathcal K_1, \rho_+)$ is again of dividing type and, in particular, $\overline{B}= - B$.
For $\mathrm{(iii)}$, $(\mathcal K_1, \rho_+)$ is of nondividing type and $U_0^{(2)} \in \BR$ follows 
from \cite[(C.37), (C.39), (C.33)]{GH}. Hence the same arguments as before yield \eqref{6.32} and \eqref{6.30}. 
\end{proof}

\begin{Rk} Given $E_j$ as in $\mathrm{(i)}$--$\mathrm{(iii)}$, we do not know if we get a solution 
	of $\snNLS_+(u)=0$ or $\snNLS_-(u)=0$ by the constraints on $A$. This has to be determined a priori,
	that is, there should be a correspondence between the location of the $E_j$'s and the defocusing/focusing nNLS equation. 
	For comparison, the stationary nonlinear Schr\"odinger potentials correspond to condition~$\mathrm{(i)}$ in the 
	defocusing case and to $\mathrm{(ii)}$ in the focusing case (see \cite[Lemmas 3.15, 3.18]{GH}). 
	\end{Rk}

\vspace{5mm}

\noindent{\bf Acknowledgments.}
 {This research   was supported by the
Austrian Science Fund (FWF) under Grants No. P29177 and V120.}

\appendix

\section{Algebro-geometric AKNS solutions} \label{App A}

Following \cite{GH}, we give a brief introduction to algebro-geometric AKNS solutions and their underlying Riemann surface
and describe the theta function representation of such solutions which we need for Lemma~\ref{L6.2}. 
The analog of these formulas for the nonlinear Schr\"odinger equation was first published by Its and Kotlyarov \cite{IK} in 1976.
Since then, many authors presented slightly varying approaches to algebro-geometric solutions of the nonlinear Schr\"odinger
and AKNS equations, see for instance Belokolos and Enol'skii \cite{BE}, Gesztesy and Ratnaseelan \cite{GR}, or Previato \cite{P}.

The stationary AKNS system \eqref{6.3} is equivalent to the stationary zero-curvature equation
\begin{equation} \label{A.13}
	\sAKNS_1(v_1,v_2) = 0 \quad \Leftrightarrow \quad - V_{2,x} + [U,V_2] = 0,
\end{equation}
where
\begin{equation}
U(z) = \begin{bmatrix} -\I z & v_1 \\
		v_2 & \I z \end{bmatrix}, \quad 
V_2(z) = \I \begin{bmatrix} - G_2(z) & F_1(z) \\
			- H_1(z) & G_2(z) \end{bmatrix}.  
\end{equation}
The polynomials $G_2(z)$, $F_1(z)$, and $H_1(z)$ incorporate the constants $c_{\ell}$,
\begin{align}\label{6.6}
	G_2(z) & = z^2+\frac{1}{2} v_1 v_2 + c_1 z + c_2, \\ \label{6.7}
	F_1(z) & = - \I v_1z + \frac{1}{2} v_{1,x}+c_1 (-\I v_1), \\ \label{6.8}
	H_1(z) & = \I v_2z + \frac{1}{2} v_{2,x}+c_1 (\I v_2).
\end{align}
The stationary zero-curvature equation in \eqref{A.13} yields that 
\begin{equation}
	\big(G_2^2 - F_1H_1\big)_x =0, 
\end{equation}
and hence $G_2^2 - F_1H_1$ is $x$-independent, implying
$G_2^2 - F_1H_1 = R_4$,
where the integration constant $R_4$ is a monic polynomial of degree $4$. If $E_0, \dots, E_3$ denote its zeros, then
\begin{equation}
R_4(z)= \prod_{m=0}^3(z-E_m), \quad \{E_m\}_{m=0}^{3} \subset \BC.	
	\end{equation}	
In this manner we can associate a hyperelliptic curve $\mathcal K_1$ of genus $1$ with \eqref{6.3} defined by
\begin{equation} \label{A.20}
	\mathcal K_1: \mathcal F_1(z,y)= y^2 - R_4(z)=0.
\end{equation}
The curve $\mathcal K_1$ is compactified by joining two points at infinity, $P_{\infty_\pm}$, $P_{\infty_+}\neq P_{\infty_-}$;
we denote the compactification again by $\mathcal K_1$. Points $P$ on $\mathcal K_1 \setminus \{P_{\infty_+}, P_{\infty_-}\}$
are represented as pairs $P=(z,y)$, where $y(\cdot)$ is the meromorphic function on $\mathcal K_1$ satisfying $\mathcal F_1(z,y)=0$.
The complex structure on $\mathcal K_1$ is then defined in the usual way (see for example \cite[App. C]{GH}). Hence 
$\mathcal K_1$ becomes a two-sheeted hyperelliptic Riemann surface of genus $1$.
We emphasize that by fixing the curve $\mathcal K_1$ (i.e., by fixing $E_0, \dots, E_3$), the integration constants 
$c_1$, $c_2$ in \eqref{6.3} are uniquely determined,
\begin{equation} \label{A.18}
c_1 = - \frac 12 (E_0 + \dots + E_3), \qquad c_2 = - \frac {c_1^2}{8} + \frac 12 \sum_{m,n=0; m<n}^{3} E_m E_n.  
\end{equation} 
Let $\mu(x)$ and $\nu(x)$ denote the zeros of $F_1(z)$ and $H_1(z)$ in \eqref{6.7} and \eqref{6.8},
\begin{equation}
	F_1(z) = - \I v_1 (z - \mu), \qquad H_1(z) = \I v_2 (z - \nu).
\end{equation}	
We lift $\mu(x)$ and $\nu(x)$ to $\mathcal K_1$ by defining
$$
\hat \mu(x) = (\mu(x), G_2(\mu(x),x)) \in \mathcal K_1, \qquad \hat \nu(x) = (\nu(x), -G_2(\nu(x),x)) \in \mathcal K_1.
$$
Choose a homology basis $\{a_1, b_1\}$ on $\mathcal K_1$ and denote by $\omega_1$ the corresponding normalized holomorphic 
differential, that is,
\begin{equation}
	\int_{a_1} \omega_1 = 1, \quad \int_{b_1} \omega_1 = \tau \in \BC.
\end{equation}
Note that $\im(\tau)>0$. 
Let $\Xi = \frac \tau 2 + \frac 12$ be the Riemann constant. 
The Riemann theta function associated with $\mathcal K_1$ is given by
\begin{equation}
	\theta (z) = \sum_{m \in \mathbb Z} \exp(2 \pi \I m z + \pi \I m^2 \tau).
\end{equation}
Without loss of generality we choose the branch point $P_0=(E_0,0)$ as a base point. 
Let $\omega_{P_{\infty_\pm,0}}^{(2)}$ be a normalized differential of the second kind satisfying 
\begin{equation}
	\int_{a_1} \omega_{P_{\infty_\pm,0}}^{(2)} = 0, \qquad 
	\omega_{P_{\infty_\pm,0}}^{(2)} = (\zeta^{-2} + O(1))d\zeta\ \mbox{as $P \to P_{\infty_\pm}$},
\end{equation}
where $\zeta$ denotes the local coordinate $\zeta=1/z$ for $P$ near $P_{\infty_\pm}$. 
Then 
\begin{equation} \label{A.23}
\int_{P_0}^P \left(\omega_{P_{\infty_+,0}}^{(2)}- \omega_{P_{\infty_-,0}}^{(2)}\right) = 
\mp \left(\zeta^{-1} + \frac{e_0}{2} + e_1 \zeta + O(\zeta^2)\right)\ \mbox{as $P \to P_{\infty_\pm}$}.
\end{equation}
In addition, we denote the $b_1$-period of this difference by
\begin{equation}
U_0^{(2)} = \frac{1}{2 \pi \I} \int_{b_1} \left(\omega_{P_{\infty_+,0}}^{(2)}- \omega_{P_{\infty_-,0}}^{(2)}\right).
\end{equation}
Finally, we turn to divisors, the Jacobi variety, and the Abel map for divisors in our setting. 
A divisor $\mathcal D$ on $\mathcal K_1$
is a map $\mathcal D : \mathcal K_1 \to \mathbb Z$, where $\mathcal D(P)\neq 0$ for only finitely many $P \in \mathcal K_1$.
We define the positive divisor $\mathcal D_Q$ by
$$
\mathcal D_Q : \mathcal K_1 \to \mathbb N_0, \quad P \mapsto \mathcal D_Q(P) = \left\{ \begin{array}{cc} 1 & \mbox{if } P=Q, \\
0 & \mbox{if } P \neq Q, 
\end{array} \right.
\quad Q \in \mathcal K_1,
$$
and denote the set of all divisors on $\mathcal K_1$ by $\textrm{Div}(\mathcal K_1)$.
The Jacobi variety $J(\mathcal K_1)$ of $\mathcal K_1$ is defined by $J(\mathcal K_1)=\BC / L$, where
$L$ is the period lattice $L= \{z \in \BC \mid z = n + m + \tau,\ n,m \in \mathbb Z\}$. The Abel map for divisors 
is then defined by
\begin{equation}
\alpha_{P_0}: \textrm{Div}(\mathcal K_1) \to J(\mathcal K_1), \quad \mathcal D \mapsto \alpha_{P_0}(\mathcal D) = 
\sum_{P \in K_1} \mathcal D (P) \int_{P_0}^{P} \omega_1.
\end{equation}
With these quantities at hand, the algebro-geometric AKNS solutions admit the following representation in terms of
Riemann theta functions, compare \cite[Thm. 3.11]{GH}.
\begin{Tm}\label{ThA.1}
	Suppose that $v_1, v_2 \in C^\infty(\Omega)$ are nonzero and satisfy the stationary AKNS system 
	\eqref{6.3} on $\Omega$. In addition, assume the affine part of $\mathcal K_1$ to be nonsingular and 
	let $x, x_0 \in \Omega$, where $\Omega \subseteq \BR$ is an open interval. Then
\begin{align} \label{6.17}
v_1(x) &= C_1 \frac {\theta (A + Bx- \Delta)}{\theta(A+Bx)}\exp(\I e_0 x), \\ \label{6.18}
	v_2(x) &= C_2 \frac {\theta (A + Bx+ \Delta)}{\theta(A+Bx)}\exp(-\I e_0 x),    
\end{align}
where
\begin{align} \label{6.19}
	A & = \Xi - \int_{P_0}^{P_{\infty_+}} \omega_1 + \I U_0^{(2)} x_0 + \alpha_{P_0}(\mathcal D_{\hat \mu(x_0)}),\\ \label{A.28}
	B & = - \I U_0^{(2)}, \qquad \Delta = \int_{P_{\infty_+}}^{P_{\infty_-}} \omega_1. 
\end{align}	  
The constants $e_0 \in \BC$ and $\Delta, B$ are uniquely determined by $\mathcal K_1$ (and its homology basis), the constant
$A$ is in one-to-one correspondence with the Dirichlet datum $\hat \mu(x_0)$ at the point $x_0$. 
The constants $C_1, C_2 \in \BC$ are given by
\begin{align} \label{6.21}
	C_1 &= v_1(x_0) \frac {\theta (A + Bx_0)}{\theta(A+Bx_0 - \Delta)}\exp(-\I e_0 x_0),     \\
	C_2 &= \frac{4}{v_1(x_0)\omega_0^2} \frac{\theta(A+Bx_0 - \Delta)}{\theta (A + Bx_0)}\exp(\I e_0 x_0),
\end{align}
and satisfy the constraint
\begin{equation}
C_1 C_2 = \frac{4}{\omega_0^2}.	
\end{equation}
\end{Tm}
Note that the free constant $v_1(x_0)$ in \eqref{6.21} cannot be determined since the AKNS equations are invariant 
with respect to scale transformations, $(v_1(x,t), v_2(x,t)) \mapsto (a v_1(x,t), a^{-1}v_2(x,t))$ for 
$a \in \BC \setminus \{0\}$.

We conclude this appendix with the following result used in the characterization of 
algebro-geometric nonlocal NLS solutions. The genus $g=1$ case of \cite[Theorem A.36 (i)]{GH} reads
\begin{Tm}\label{ThA.2}
	Let $(\mathcal K_1, \rho)$ be a symmetric Riemann surface, i.e., let $\rho$ be an antiholomorphic involution on $\mathcal K_1$. There exists a canonical homology basis $\{a_1, b_1\}$ on 
	$\mathcal K_1$ with intersection index $a_1 \circ b_1 = 1$ such that 
	the $2 \times 2$ matrix $S$ of complex conjugation of the action of $\rho$ on $H_1(\mathcal K_1, \mathbb Z)$ in this 
	basis is given by 	$$	S=\begin{pmatrix}1 & 0 \\ 0 & -1\end{pmatrix},	$$	that is, 	
	$$
	(\rho(a_1), \rho(b_1)) = (a_1,b_1) \begin{pmatrix} 1 & 0 \\ 0 & -1 \end{pmatrix} =
	(a_1, - b_1).
	$$
\end{Tm}



\begin{flushright}

J. Michor, \\
Fakult\"at f\"ur Mathematik, Universit\"at Wien, \\
Oskar-Morgenstern-Platz 1, A-1090 Vienna, Austria\\
e-mail: Johanna.Michor@univie.ac.at

\vspace{0.5em}

A.L. Sakhnovich,\\
Fakult\"at f\"ur Mathematik, Universit\"at Wien, \\
Oskar-Morgenstern-Platz 1, A-1090 Vienna, Austria\\
e-mail: oleksandr.sakhnovych@univie.ac.at
\end{flushright}

\end{document}